\documentclass{amsart}
\usepackage{graphicx}
\usepackage{amssymb}
\usepackage{amsfonts}
\swapnumbers
\newtheorem{theorem}{Theorem}[section]
\newtheorem{lemma}[theorem]{Lemma}
\newtheorem{corollary}[theorem]{Corollary}
\newtheorem{proposition}[theorem]{Proposition}
\theoremstyle{definition}

\newtheorem{assumption}[theorem]{Assumption}
\newtheorem{remark}[theorem]{Remark}

\numberwithin{equation}{section}
 \theoremstyle{plain}    
 
 \numberwithin{equation}{section} 
 \numberwithin{figure}{section} 
 \theoremstyle{plain}    
 \theoremstyle{plain}    
 \theoremstyle{remark}    
 \newtheorem*{acknowledgement*}{Acknowledgement} 

\newcommand{\cA}{{\mathcal A}}
\newcommand{\cB}{{\mathcal B}}

\newcommand{\cF}{{\mathcal F}}
\newcommand{\cG}{{\mathcal G}}
\newcommand{\cH}{{\mathcal H}}

\newcommand{\cN}{{\mathcal N}}

\newcommand{\te}{{\theta}}

\newcommand{\Om}{{\Omega}}
\newcommand{\om}{{\omega}}
\newcommand{\ve}{{\varepsilon}}
\newcommand{\del}{{\delta}}

\newcommand{\gam}{{\gamma}}
\newcommand{\Gam}{{\Gamma}}
\newcommand{\vf}{{\varphi}}
\newcommand{\vp}{{\varpi}}
\newcommand{\io}{{\iota}}
\newcommand{\up}{{\upsilon}}

\newcommand{\sig}{{\sigma}}
\newcommand{\al}{{\alpha}}
\newcommand{\be}{{\beta}}
\newcommand{\ka}{{\kappa}}

\newcommand{\bbN}{{\mathbb N}}

\newcommand{\bbR}{{\mathbb R}}

\newcommand{\bbI}{{\mathbb I}}

\begin{document}
\title[]{A nonconventional strong law of large numbers\\
 and fractal dimensions of some multiple recurrence sets}%
 \vskip 0.1cm 
 \author{ Yuri Kifer\\
 \vskip 0.1cm
Institute of Mathematics\\
The Hebrew University of Jerusalem}%
\email{kifer\@@math.huji.ac.il}
\address{Institute of Mathematics, Hebrew University, Jerusalem 91904,\linebreak
 Israel}


\thanks{ }
\dedicatory{  }
 \date{\today}
\begin{abstract}\noindent
We provide conditions which yield a strong law of large numbers for
expressions of the form
$1/N\sum_{n=1}^{N}F\big(X(q_1(n)),\cdots, X(q_\ell(n))\big)$
where $X(n),n\geq 0$'s is a sufficiently fast 
mixing vector process with some moment conditions and stationarity properties,
$F$ is a continuous function with polinomial growth and certain regularity 
properties and $q_i,i>m$ are positive functions taking on integer values on
integers with some growth conditions. Applying these results
 we study certain multifractal formalism type questions concerning
Hausdorff dimensions of some sets of numbers with prescribed asymptotic
frequencies of combinations of digits at places $q_1(n),...,q_\ell(n)$.  

\end{abstract}
\maketitle
\markboth{Yu.Kifer}{Strong LLN} 
\renewcommand{\theequation}{\arabic{section}.\arabic{equation}}
\pagenumbering{arabic}

\section{Introduction}\label{sec1}\setcounter{equation}{0}

Nonconventional ergodic theorems which attracted substantial attention 
in ergodic theory (see, for instance, \cite{Be}, \cite{Fu} and \cite{As}) 
studied the limits
 of expressions having the form
$1/N\sum_{n=1}^NT^{q_1(n)}f_1\cdots
 T^{q_\ell (n)}f_\ell$ where $T$ is a
weakly mixing measure
 preserving transformation, $f_i$'s are bounded 
measurable functions
 and $q_i$'s are polynomials taking on integer values on
the
 integers. While, for instance, \cite{Be} and \cite{Fu} were interested
in $L^2$ convergence, other papers such as \cite{As} provided conditions
for almost sure convergence in such ergodic theorems.
Originally, these results were motivated by applications to 
multiple recurrence for dynamical systems taking functions $f_i$ being 
indicators of some measurable sets. 

Introducing stronger mixing or weak dependence conditions enabled us in
\cite{Ki3} and \cite{KV} to obtain central limit theorems and invariance 
principles for even more general expressions of the form
\begin{equation}\label{1.1}
\frac 1{\sqrt N}\sum_{n=1}^{[Nt]}\big( F(X(q_1(n)),...,X(q_\ell(n))-\bar F\big)
\end{equation}
where $X(n),\, n\geq 0$ is a sufficiently fast mixing vector valued process
with some moment conditions and stationarity properties, $F$ is a locally 
H\" older continuous function with polinomial growth, $\bar F=\int Fd(\mu\times
\cdots\times\mu)$ and $\mu$ is the distribution of $X(0)$. In order to ensure
existence of limiting variances and covariances we had to impose another 
assumption concerning the functions $q_j(n),\, j\geq 1$ saying that $q_j(n)=jn$
for $j=1,...,k$ while $q_j(n),\, j\geq k$ are positive functions taking on
integer values on integers with some (faster than linear) growth conditions.

In this paper we are concerned with strong laws of large numbers (SLLN) for 
expressions of the form
\begin{equation}\label{1.2}
\frac 1N\sum_{n=1}^N F(X(q_1(n)),...,X(q_\ell(n))
\end{equation}
which can be proved under milder conditions that those required for central
limit theorem type results. We still impose some mixing or weak dependence
conditions but now the functions $q_j(n),\, n\geq 1$ are allowed to be of much
more general form than in \cite{KV}, in particular, because we do not have to
take care about limiting variances. Recall, that the machinery of 
nonconventional ergodic theorems employed in \cite{Be}, \cite{Fu}, \cite{As}
and other
papers can only work when the functions $q_j,\,j=1,...,\ell$ are polinomials
while our methods do not require any algebraic structure of them. We pay a
price for this, namely, imposing stronger mixing assumptions which are satisfied
though for important classes of stochastic processes and dynamical systems.

In order to obtain our strong laws of large numbers we represent the sum in 
(\ref{1.2}) as a sum of certain mixingales and then rely on the SLLN for 
mixingales obtained in \cite{ML}. Another approach which 
works in this situation under more or less the same assumptions is a 
martingale approximation similar to \cite{KV} together with a SLLN
 for martingales (see, for instance, Section 2.6 in \cite{HH}).

Among more specific applications of our setup we can consider $F(x_1,...,
x_\ell)=x_1^{(1)}\cdots x_\ell^{(\ell)},$ $x_j=(x_j^{(1)},...,x_j^{(\ell)}),$
$X(n)=(X_1(n),...,X_|ell(n)),$ $X_j(n)=\bbI_{A_j}(T^nx)$ for a dynamical
system $\{ T^n\}$ or $X_j(n)=\bbI_{A_j}(\xi_n)$ for a Markov chain $\{\xi_n\}$
where $\bbI_A$ is the indicator of a set $A$. Then the expression (\ref{1.2})
measures the frequency of arrivals of $T^nx$ or of $\xi_n$ to the sets $A_j$
at the respective times $q_j(n)$. Recall, that the $m$-base and continued 
fraction expansions can be obtained via the multiplication by $m$ and the
Gauss transformations, i.e. $Tx=\{ mx\}$ and $Tx=\{ 1/x\}$, respectively, 
which are both exponentially fast $\psi$-mixing with respect to many invariant
measures (see \cite{He} and \cite{Aa}) and satisfy
our assumptions. Denote by $\zeta_j(x)$ the $j$-th digit of $x\in[0,1)$ in
one of these expansions. Then we can study the frequency of $k$-th such that
the $\ell$-tuple $(\zeta_{q_1(k)}(x),...,\zeta_{q_\ell(k)}(x))$ coincides 
with a prescribed $\ell$-tuple of digits $(a_1,...,a_\ell)$. For a full
Lebesgue measure of points $x\in[0,1)$ such frequencies are determined by
our SLLN and other frequencies may occur only for $x$
belonging to sets of zero measure. This leads to an interesting question 
about Hausdorff dimensions of such exceptional sets which we study in the 
last section of this paper.

\section{Preliminaries and main results}\label{sec2}\setcounter{equation}{0}

Our setup consists of a  $\wp$-dimensional stochastic process
$\{X(n),  n=0,1,...\}$ on a probability space $(\Om,\cF,P)$ and of a family
of $\sig$-algebras \hbox{$\cF_{kl}\subset\cF,\, 0\leq k\leq l\leq\infty$} 
where we assume that $\cF_{00}$ is a trivial $\sig$-field and
 $\cF_{kl}\subset\cF_{k'l'}$ if $k'\leq k$ and
$l'\geq l$. We extend $\cF_{kl}$ also to negative $k\geq -\infty$ by
defining $\cF_{kl}=\cF_{0l}$ for $k<0$ and $l\geq 0$. The dependence between 
two sub $\sig$-algebras $\cG,\cH\subset\cF$ is measured often via the quantities
\begin{equation}\label{2.1}
\varpi_{q,p}(\cG,\cH)=\sup\{\| E\big [g|\cG\big]-E[g]\|_p:\, g\,\,\mbox{is}\,\,
\cH-\mbox{measurable and}\,\,\| g\|_q\leq 1\},
\end{equation}
where the supremum is taken over real functions and $\|\cdot\|_r$ is the
$L^r(\Om,\cF,P)$-norm. Then more familiar $\al,\rho,\phi$ and $\psi$-mixing 
(dependence) coefficients can be expressed in the form (see \cite{Bra}, Ch. 4 ),
\begin{eqnarray*}
&\al(\cG,\cH)=\frac 14\varpi_{\infty,1}(\cG,\cH),\,\,\rho(\cG,\cH)=\varpi_{2,2}
(\cG,\cH)\\
&\phi(\cG,\cH)=\frac 12\varpi_{\infty,\infty}(\cG,\cH)\,\,\mbox{and}\,\,
\psi(\cG,\cH)=\varpi_{1,\infty}(\cG,\cH).
\end{eqnarray*}
The relevant quantities in our setup are
\begin{equation}\label{2.2}
\varpi_{q,p}(n)=\sup_{k\geq 0}\varpi_{q,p}(\cF_{-\infty,k},\cF_{k+n,\infty})
\end{equation}
and accordingly
\[
\al(n)=\frac{1}{4}\varpi_{\infty,1}(n),\,\rho(n)=\varpi_{2,2}(n),\,
\phi(n)=\frac 12\varpi_{\infty,\infty}(n)\,\,\mbox{and}\,\,\ \psi(n)=
\varpi_{1,\infty}(n).
\]
Our assumptions will require certain speed of decay as $n\to\infty$ of both
the mixing rates $\varpi_{q,p}(n)$ and the approximation rates defined by
\begin{equation}\label{2.3}
\beta_p (n)=\sup_{m\geq 0}\|X(m)-E\big (X(m)|\cF_{m-n,m+n}\big)\|_p.
\end{equation}

Furthermore, we do not require stationarity of the
process $X(n), n\geq 0$ assuming only that the distribution $\mu$ of $X(n)$ 
does not depend on $n$ which we write for further references by 
\begin{equation}\label{2.4}
X(n)\stackrel {d}{\sim}\mu
\end{equation}
where $Y\stackrel {d}{\sim}Z$ means that $Y$ and $Z$ have the same distribution.

Next,let $F= F(x_1,...,x_\ell),\, x_j\in\bbR^{\wp}$ be a function on 
$\bbR^{\wp\ell}$ such that for some $\iota,K>0,\ka\in (0,1]$ and all 
$x_i,y_i\in\bbR^{\wp}, i=1,...,\ell$, 
\begin{equation}\label{2.5}
|F(x_1,...,x_\ell)-F(y_1,...,y_\ell)|\leq K\big(1+\sum^\ell_{j=1}|x_j|^\iota+
\sum^\ell_{j=1} |y_j|^\iota \big)\sum^\ell_{j=1}|x_j-y_j|^\ka
\end{equation}
and 
\begin{equation}\label{2.6} 
|F(x_1,...,x_\ell)|\leq K\big( 1+\sum^\ell_{j=1}|x_j|^{\iota} \big).
\end{equation}
Our assumptions on $F$ are motivated by the desire to include, for instance,
products $F(x_1,...,x_\ell)=x_{11}x_{22}\cdots x_{\ell\ell}$, where 
$x_i=(x_{i1},...,x_{i\ell})\in\bbR^\ell$, which are important in the study
of multiple recurrence as described in Introduction.

Our setup includes also a sequence of positive functions 
$q_1(n)< q_2(n) <\cdots < q_\ell(n)$ taking on integer values on integers and 
such that for some positive $\ve\leq 1$,
\begin{equation}\label{2.7}
q_i(n)\geq q_{i-1}(n)+\ve n,\, i=2,...,\ell\,\,\mbox{and}\,\, q_i(n+1)\geq
q_i(n)+\ve\,\,\mbox{for all}\,\, n\geq 1.
\end{equation}
In order to give a detailed statement of our main result as well as for its
proof it will be essential to represent the function $F= F(x_1,x_2,\ldots,
x_\ell)$ in the form 
\begin{equation}\label{2.8}
F=F_0+F_1(x_1)+\cdots+F_\ell(x_1, x_2,\ldots, x_\ell)
\end{equation}
where
\begin{equation}\label{2.9}
F_0=\bar F=\int F(x_1,...,x_\ell)\,d\mu(x_1)\cdots d\mu(x_\ell),
\end{equation} 
\begin{eqnarray}\label{2.10}
&F_i(x_1,\ldots, x_i)=\int F(x_1,x_2, \ldots, x_\ell)\ d\mu (x_{i+1})\cdots 
d\mu(x_\ell)\\
&\quad -\int F(x_1,x_2, \ldots, x_\ell) \,d\mu (x_i)\cdots d\mu(x_\ell)\nonumber
\end{eqnarray}
for $0<i<\ell$ and
\[
F_\ell(x_1,x_2, \ldots, x_\ell)=F(x_1,x_2, \ldots, x_\ell) -\int F(x_1,x_2, 
\ldots, x_\ell)\, d\mu(x_\ell)
\]
which ensures, in particular, that
\begin{equation}\label{2.11}
\int F_i(x_1, x_2,\ldots,x_{i-1}, x_i)\,d\mu(x_i)\equiv 0 \quad\forall 
\quad x_1, x_2,\ldots, x_{i-1}.
\end{equation}
These enable us to write
\begin{equation}\label{2.12}
S(N)=\sum_{n=1}^N F(X(q_1(n)),...,X(q_\ell(n)))=\sum_{i=0}^\ell S_i(N)
\end{equation}
where $S_0(N)=N\bar F$ and for $1\leq i\leq\ell $,
\begin{equation}\label{def2.13}
S_i(N)=\sum_{1\leq n\leq N} F_i(X(q_1(n)),\ldots, X(q_i(n))).
\end{equation}

Following \cite{ML} we say that a sequence $\{ a_n,\, n\geq 0\}$ is of size
$-1/2$ if there exists a positive eventually nondecreasing sequence 
$\{ L_n,\, n\geq 0\}$ such that 
\[
\sum_{n\geq 0}(nL_n)^{-1}<\infty,\, L_n-L_{n-1}=O(L_n/n)\,\,\mbox{and}\,\,
a_n=O\big((n^{1/2}L_n)^{-1}\big).
\]
For instance, any sequence with asymptotics $O\big(n^{1/2}\log n (\log\log
n)^{1+\del}\big)^{-1}$ for some $\del>0$ is of size $-1/2$.
For each $r>0$ set
\begin{equation}\label{2.14}
\gamma_r^r = \|X\|_r^r= E|X(n)|^r  =
\int \|x\|^r d\mu .
\end{equation}
Our main result relies on

\begin{assumption}\label{ass2.1} With $d=(\ell-1)\wp$ there exist $p,q\geq 1$
and $\theta,m>0$ such that $\theta<\ka-\frac dp,$ 
\begin{equation}\label{2.15}
\frac 12\geq\frac 1p+\frac {\iota+2}m+\frac \theta q\,\,\,\mbox{and}\,\,\,
\gam_m+\gam_{2q(\io+2)}<\infty
\end{equation}
and the sequence $\vp_{p,q}(n)+\be_q^\theta(n),\, n\geq 1$ is of size 
$-1/2$.
\end{assumption}

\begin{theorem}\label{thm2.2}
Suppose that Assumption \ref{ass2.1} holds true. Then with probability one
\begin{equation}\label{2.16}
\lim_{N\to\infty}\frac 1NS(N)=\bar F.
\end{equation}
\end{theorem}

Our method relies on estimates from \cite{KV} which enable us to view 
for each $i\geq 1$ the sequence of pairs $\{ F_i(X(q_1(n)),...,X(q_i(n))),
\,\cF_{-\infty,q_i(n)}\}_{n=1}^\infty$ as a mixingale sequence, and so
a strong law of large numbers for mixingales from \cite{ML} can be employed.
This gives an almost sure convergence of $\frac 1NS_i(N)$ to 0 and by 
(\ref{2.12}) Theorem \ref{thm2.2} follows.
Another approach which works in our situation is to rely on a martingale 
approximation of $S_i(N)$ similarly to \cite{KV} and then to employ a
strong law of large numbers for martingales (see, for instance, Section 2.6
in \cite{HH}). This method has to deal with approximations of 
$F_i(X(q_1(n)),...,X(q_i(n)))$ by their conditional expectations and in order
to avoid double limits as in \cite{KV} we can make this approximations with
increasing in $n$ precision.

In order to understand our assumptions observe that $\varpi_{q,p}$
 is non-increasing in $q$ and non-decreasing in $p$. Hence,
 for any pair $p,q\geq 1$,
 \[
 \varpi_{q,p}(n)\leq\psi(n).
 \]
 Furthermore, by the real version of the Riesz--Thorin interpolation 
 theorem (see, for instance, \cite{Ga}, Section 9.3) if 
 $\del\in[0,1],\, 1\leq p_0,p_1,q_0,q_1\leq\infty$ and
 \[
 \frac 1p=\frac {1-\del}{p_0}+\frac \del{p_1},\,\,\frac 1q=\frac
 {1-\del}{q_0}+\frac \del{q_1}
 \]
 then
 \begin{equation*}
\varpi_{q,p}(n)\le 2(\varpi_{q_0,p_0}(n))^{1-\del}
(\varpi_{q_1,p_1}(n))^\del.
\end{equation*}
Since, clearly, $\varpi_{q_1,p_1}\leq 2$ for any $q_1\geq p_1$ it follows
for pairs $(\infty,1)$, $(2,2)$ and $(\infty,\infty)$ that for all 
$q\geq p\geq 1$,
 \begin{eqnarray*}
&\varpi_{q,p}(n)\le (2\alpha(n))^{\frac{1}{p}-\frac{1}{q}},\,
 \varpi_{q,p}(n)\le 2^{1+\frac 1p-\frac 1q}(\rho(n))^{1-\frac 1p+\frac 1q}\\
&\mbox{and}\,\,\varpi_{q,p}(n)\le 2^{1+\frac 1p}(\phi(n))^{1-\frac 1p}.
\end{eqnarray*}
We observe also that by the H\" older inequality for $q\geq p\geq 1$
and $\alpha\in(0,p/q)$,
\begin{equation*}
\beta(q,r)\le 2^{1-\alpha}  [\beta(p,r)]^\alpha \gamma^{1-\al}_{\frac{pq(1-\al)}
{p-q\al}}
\end{equation*}
with $\gamma_r$ defined in (\ref{2.14}). Thus, we can formulate 
Assumption \ref{ass2.1} in terms of more familiar $\alpha,\,\rho,\,\phi,$
and $\psi$--mixing coefficients and with various moment conditions. 

The conditions of Theorem \ref{thm2.2} hold true for 
many important models. Let, for instance, $\xi_n$ be a Markov chain on a
space $M$ satisfying the Doeblin condition (see, for instance,
\cite{IL}, p.p. 367--368) and $f_j,\, j=1,...,\ell$ be bounded measurable
functions on the space of sequences $x=(x_i,\, i=0,1,2,...,\, x_i\in M)$ such
that $|f_j(x)-f_j(y)|\leq Ce^{-cn}$ provided $x=(x_i),\, y=(y_i)$
and $x_i=y_i$ for all $i=0,1,...,n$ where $c,C>0$ do not depend on
$n$ and $j$. In fact, some polinomial decay in $n$ will suffice here, as
well. Let $X(n)=(X_1(n),...,X_\ell(n))$ with $X_j(n)=
f_j(\xi_n,\xi_{n+1},\xi_{n+2},...)$ and take $\sig$-algebras $\cF_{kl},\,
 k<l$ generated by $\xi_k,\xi_{k+1},...,\xi_l$ then our condition 
will be satisfied considering $\{\xi_n,\, n\geq 0\}$ with its invariant
measure as a stationary process. In fact, our conditions hold true for a
more general class of processes, in particular, for Markov chains whose
transition probability has a spectral gap which leads to an exponentially 
fast decay of the $\rho$-mixing coefficient.

Important classes of processes satisfying our conditions come from
dynamical systems. Let $T$ be a $C^2$ Axiom A diffeomorphism (in
particular, Anosov) in a neighborhood of an attractor or let $T$ be
an expanding $C^2$ endomorphism of a compact Riemannian manifold $M$ (see
\cite{Bo}), $f_j$'s be H\" older continuous functions and let
$X(n)=(X_1(n),...,X_\ell(n))$ with $X_j(n)=f_j(T^nx)$. Here the probability
space is $(M,\cB,\mu)$ where $\mu$ is a Gibbs invariant measure corresponding
to some H\"older continuous function and $\cB$ is the Borel $\sig$-field. Let
$\zeta$ be a finite Markov partition
for $T$ then we can take $\cF_{kl}$ to be the finite $\sig$-algebra
generated by the partition $\cap_{i=k}^lT^i\zeta$. In fact, we can
take here not only H\" older continuous $f_j$'s but also indicators
of sets from $\cF_{kl}$. A related example corresponds to $T$ being
a topologically mixing subshift of finite type which means that $T$
is the left shift on a subspace $\Xi$ of the space of one-sided
sequences $\xi=(\xi_i,i\geq 0), \xi_i=1,...,l_0$ such that $\xi\in\Xi$
if $\pi_{\xi_i\xi_{i+1}}=1$ for all $i\geq 0$ where $\Pi=(\pi_{ij})$ 
is an $l_0\times l_0$
matrix with $0$ and $1$ entries and such that $\Pi^n$ for some $n$
is a matrix with positive entries. Again, we have to take in this
case $f_j$ to be H\" older continuous bounded functions on the
sequence space above, $\mu$ to be a Gibbs invariant measure
corresponding to some H\" older continuous function and to define
$\cF_{kl}$ as the finite $\sig$-algebra generated by cylinder sets
with fixed coordinates having numbers from $k$ to $l$. The
exponentially fast $\psi$-mixing is well known in
the above cases (see \cite{Bo}). Among other dynamical systems with
exponentially fast $\psi$-mixing we can mention also the Gauss map
$Tx=\{1/x\}$ (where $\{\cdot\}$ denotes the fractional part) of the
 unit interval with respect to the Gauss measure $G(\Gam)=\frac 1{\ln 2}
 \int_\Gam\frac 1{1+x}dx$ (see \cite{He}), as well as with respect to
 many other Gibbs invariant measures (see \cite{Aa}).
 The latter enables us to consider the number $N_a(x,n)$,
$a=(a_1,...,a_\ell)$
of $k$'s between 0 and $n$ such that the $q_j(k)$-th digit of the
continued fraction of $x$ equals certain integer $a_j,j=1,...,\ell$. Then
Theorem \ref{thm2.2} implies a strong law of large numbers for $N_a(x,n)$ 
considered
as a random variable on the probability space $((0,1],\cB, G)$. In fact, our
results rely only on sufficiently fast
$\al$ or $\rho$-mixing which holds true for wider classes of dynamical system,
in particular, those with a spectral gap (such as many one dimensional
not necessarily uniformly expanding maps) which ensures an exponentially fast
$\rho$-mixing. We will show how to derive from Theorem \ref{thm2.2} the 
following result.

\begin{corollary}\label{cor2.3}
Let $T$ be either a $C^2$ Axiom A diffeomorphism on a compact Riemannian
manifold $M$ considered in a neighborhood of an attractor or a $C^2$ expanding
endomorphisms of a compact Riemannian manifold $M$ or the Gauss map of the unit 
interval and let $\mu$ be an equilibrium state (Gibbs measure) corresponding to
a H\" older continuous function in the first two cases or an exponentially fast 
$\psi$-mixing $T$-invariant (in particular, Gauss') measure (see Corollary 
4.7.8 in \cite{Aa}) in the latter case. 
Let $X_j(n)=f_j(T^nx),\, j=1,...,\ell$ where $f_j$ is 
either a continuous function or $f_j(x)=\bbI_{\Gam_j}(x)$ where $\Gam_j$
is a measurable set whose boundary $\partial\Gam_j$ has zero $\mu$-measure.
Finally, let $F=F(x_1,...,x_\ell)$ satisfies conditions of Theorem \ref{thm2.2}
which means just that $F$ is H\" older continuous since its arguments are
bounded here. Then the conclusion of Theorem \ref{thm2.2} holds true.
\end{corollary}

Next, we discuss a continuous time version of our theorem. 
Our continuous time setup consists of a $\wp$-dimensional process
$X(t),\, t\geq 0$ on a probability space $(\Om,\cF,P)$ whose one
dimensional distributions do not depend on time and of a family of 
$\sig$-algebras $\cF_{st}\subset\cF,\,-\infty\leq s\leq t\leq\infty$ such 
that $\cF_{st}\subset\cF_{s't'}$ if $s'\leq s$ and $t'\geq t$. For all 
$t\geq 0$ we set
\begin{equation}\label{2.17}
\varpi_{q,p}(t)=\sup_{s\geq 0}\varpi_{q,p}(\cF_{-\infty,s},\cF_{s+t,\infty})
\end{equation}
and
\begin{equation}\label{2.18}
\beta (p,t)=\sup_{s\geq 0}\|X(s)-E\big [X(s)|\cF_{s-t,s+t}\big]\|_p.
\end{equation}
where $\varpi_{q,p}(\cG,\cH)$ is defined by (\ref{2.1}). It will suffice
for our purposes to rely on Assumtion \ref{ass2.1} concerning
$\varpi_{q,p}(t)$ and $\beta (p,t)$ considered only for integer $t$. Let
$q_1(t)< q_2(t) <\cdots < q_\ell(t)$ be increasing positive functions
satisfying the conditions (\ref{2.7}) with $t$ in place of $n$. Set 
\begin{equation}\label{2.19}
S(t)=\int_0^tF(X(q_1(s)),...,X(q_\ell(s)))ds=\sum_{i=0}^\ell S_i(t)
\end{equation}
where $S_0(t)=t\bar F$,
\begin{equation}\label{def2.20}
S_i(t)=\int_0^t F_i(X(q_1(s)),\ldots, X(q_i(s)))ds
\end{equation}
and $F,\, \bar F,\, F_i$ are the same as in (\ref{2.5}), (\ref{2.6}) and
(\ref{2.8})--(\ref{2.11}). Then we obtain

\begin{corollary}\label{cor2.4}
Under the conditions above with probability one
\begin{equation*}
\lim_{t\to\infty}\frac 1tS(t)=\bar F.
\end{equation*}
\end{corollary}

Next, we discuss the fractal dimensions part of this paper. Recall that the
multifractal formalism deals with computations of Hausdorff dimensions of sets
having the form
\[
\{ x:\,\lim_{n\to\infty}\frac 1N\sum_{n=1}^Nf(T^nx)=\rho\}.
\]
In our setup it is natural to study Hausdorff dimensions of more general sets
\[
G_\rho=\{ x:\,\lim_{N\to\infty}\frac 1N\sum_{n=1}^NF(f_1(T^{q_1(n)}x),...,
f_\ell(T^{q_\ell(n)}x))=\rho\},
\]
say, under the conditions of Corollary \ref{cor2.3}. When
\[
\rho=\int...\int F(f_1(x_1),...,f_\ell(x_\ell))d\mu(x_1)\cdots d\mu(x_\ell)
\]
then $\mu(G_\rho)=1$ by Corollary \ref{cor2.3} while otherwise 
$\mu(G_\rho)=0$ and it is natural to inquire about the Hausdorff dimension
 of $G_\rho$.
 
 We will not study here this general problem but consider a more specific
 question about Hausdorff dimensions of sets of numbers with prescribed
 frequencies of specific combinations of digits in $m$-expansions. Namely,
 for any $x\in[0,1]$ and an integer $m>1$ we can write
 \[
 x=\sum_{i=1}^\infty\frac {a_{i-1}(x)}{m^i}\,\,\mbox{where}\,\, a_j(x)\in
 \{ 0,1,...,m-1\},\, j=0,1,...
 \]
 and we allow zero tails of expansions but not tails consisting of all
 $(m-1)$'s. This convention affects only a countable number of points, and
 so it does not influence Hausdorff dimensions computations. For each
 $x\in[0,1]$ and an $\ell$-word $\al=(\al_1,\al_2,...,\al_\ell)\in\{ 0,1,...,
 m-1\}^\ell$ define 
 \begin{equation}\label{2.21}
 N_\al(x,n)=\#\{ k>0,k\leq n:\, (a_{q_1(k)}(x),...,a_{q_\ell(k)}(x))=\al\}
 \end{equation}
 where $\#\Gam$ denotes the number of elements in the set $\Gam$. Denote by
 $\cA_\ell=\{ 0,1,...,m-1\}^\ell$ the set of all $\ell$-words and let
 $p_\al\geq 0,\,\al\in\cA_\ell$ satisfy $\sum_{\al\in\cA_\ell}p_\al=1$.
 For such a probability vector $p=(p_\al,\,\al\in\cA_\ell)\in\bbR^{m^\ell}$
 define
 \begin{equation}\label{2.22}
 U_p=\{ x\in(0,1):\,\lim_{n\to\infty}\frac 1nN_\al(x,n)=p_\al\,\,\mbox{for
 all}\,\,\al\in\cA_\ell\}.
 \end{equation}
 We want to deal with the question of computation of the Hausdorff dimension
 $HD(U_p)$ of $U_p$. When $\ell=1$ and $q_1(k)=k$ we arrive at the classical
 question studied in \cite{Be} and \cite{Eg} by combinatorial means and in
 \cite{Bi2} via the ergodic theory.
 
 In order to relate the limit of $n^{-1}N_\al(x,n)$ to the nonconventional
 strong law of large numbers (ergodic theorem) discussed before define
 the transformation $Tx=\{ mx\}$ where $\{\cdot\}$ denotes the fractional 
 part. Identifying 0 and 1 we can view $T$ as an expanding map of the circle.
 Now $a_i(x)=a_0(T^ix)$ and if $\al=(\al_1,\al_2,...,\al_\ell)\in\cA_\ell$
 and $\Gam_j=\{ x:\, a_0(x)=j\}$ then
 \begin{equation}\label{2.23}
 N_\al(x,n)=\sum_{k=1}^n\bbI_{\Gam_{\al_1}}(T^{q_1(k)}x)
 \bbI_{\Gam_{\al_2}}(T^{q_2(k)}x)\cdots\bbI_{\Gam_{\al_\ell}}(T^{q_\ell(k)}x).
 \end{equation}
 Taking into account that $\{\Gam_j,\, j=0,1,...,m-1\}$ is the Markov partition
 for $T$ in this simple situation we arrive at the setup of Corollary 
 \ref{cor2.3} with $F(x_1,...,x_\ell)=x_1x_2\cdots x_\ell$ and $f_j(x)=
 \bbI_{\Gam_{\al_j}}(x),\, j=1,...,\ell$. Observe that in place of the 
 dynamical systems setup described above we could rely in this situation 
 on the fact that that the digits $a_n,\, n\geq 0$ are independent identically
 distributed (i.i.d.) random 
 variables with respect to the Lebesgue measure on $[0,1]$, and so
 $\bbI_{\Gam_{\al_j}}\circ T^n=\bbI_{a_n=\al_j},\, i=1,...,\ell,n=0,1,...$
 are also i.i.d. random variables so that mixing conditions of Assumption
 \ref{ass2.1} trivially hold true. The following result answers our
 question in a specific situation. 
 
 \begin{proposition}\label{prop2.5} Suppose that $q_1(k)=k$ for all $k$ and 
 there exists a probability vector 
 $r=(r_0,r_1,...,r_{m-1})$ such that $p_\al=\prod_{i=1}^\ell r_{\al_i}$ for
 any $\al=(\al_1,...,\al_\ell)\in\cA_\ell$. Then
 \begin{equation}\label{2.24}
 HD(U_p)=\frac {-\sum_{j=0}^{m-1}r_j\ln r_j}{\ln m}
 \end{equation}
 with the convention $0\ln 0=0$.
 \end{proposition}
 
 \begin{remark}\label{rem2.6}
 In view of (\ref{2.23}) for any $T$-invariant probability measure $\mu$
 on $[0,1]$ with mixing properties fulfilling conditions of Theorem 
 \ref{thm2.2} it follows that $\mu$-almost everywhere
 \[
 \lim_{n\to\infty}\frac 1nN_\al(x,n)=\prod_{i=1}^\ell\mu(\Gam_{\al_i}).
 \]
 Hence, if $p=(p_\al,\,\al\in\cA_\ell)$ and there exists no probability
 vector $r=(r_0,r_1,...,r_{m-1})$ such that $p_\al=\prod_{i=1}^\ell 
 r_{\al_i}$ then $\mu(U_p)=0$ for any $\mu$ as above, and so such $\mu$
 cannot be used for computation of the Hausdorff dimension of $U_p$
 (by one of methods where measures are involved) which complicates the
 study in this case.
 \end{remark} 

 Now, consider a bit more complex situation. For each $x\in[0,1]$ and 
 $\al=(\al_1,\al_2,...,\al_\ell),\,\be=(\be_1,\al_2,...,\be_\ell)\in
 \{ 0,1,..., m-1\}^\ell$ set
 \begin{eqnarray*}
 &N_{\al,\be}(x,n)=\#\{ k>0,k\leq n:\, (a_{q_1(k)}(x),...,a_{q_\ell(k)}(x))
 =\al\\
 &\mbox{and}\,\,(a_{q_1(k)+1}(x),...,a_{q_\ell(k)+1}(x))=\be\}
 \end{eqnarray*}
 and for each nonnegative matrix $P=(p_{\al\be},\,\al,\be\in\cA_\ell)$
 with $\sum_{\al,\be}p_{\al\be}=1$ define
 \begin{equation}\label{2.25}
 U_P=\{ x\in(0,1):\,\lim_{n\to\infty}\frac 1nN_{\al,\be}(x,n)=p_{\al,\be}\,\,
 \mbox{for all}\,\,\al,\be\in\cA_\ell\}.
 \end{equation}
 Again, we can write $N_{\al\be}$ in the form suitable for application of
 Theorem \ref{thm2.2}, namely,
 \begin{equation}\label{2.26}
 N_{\al,\be}(x,n)=\sum_{k=1}^n\bbI_{\Gam_{\al_1\be_1}}(T^{q_1(k)}x)
 \bbI_{\Gam_{\al_2\be_2}}(T^{q_2(k)}x)\cdots\bbI_{\Gam_{\al_\ell\be_\ell}}
 (T^{q_\ell(k)}x)
 \end{equation}
where $\Gam_{ij}=\{ x:\, a_0(x)=i,\, a_1(x)=j\}$. Then we obtain the following
result.
\begin{proposition}\label{prop2.7}
Suppose that $q_1(k)=k$ and there exists a nonnegative matrix $R=(r_{ij};\,
i,j=0,1,...,m-1)$ satisfying the following conditions:

 (i) some power of $R$ is a positive matrix;
(ii) $\sum_{i,j}r_{ij}=1$, $q_i=\sum_{j=0}^{m-1}r_{ij}=
\sum_{j=0}^{m-1}r_{ij}$;
(iii) $p_{\al\be}=\prod_{i=1}^\ell r_{\al_i\be_i}$.

Then $q=(q_0,q_1,...,q_{m-1})$ is a positive stationary vector of the 
$m\times m$ irredicible aperiodic probability matrix $Q=(q_{ij})$, 
$q_{ij}=q_i^{-1}r_{ij}$ and under the convention $0\ln 0=0$,
\begin{equation}\label{2.27}
HD(U_P)=\frac {-\sum_{i,j=0}^{m-1}q_iq_{ij}\ln q_{ij}}{\ln m}.
\end{equation}
\end{proposition}

\begin{remark}\label{rem2.8}
Somewhat surprisingly Proposition \ref{prop2.5} and \ref{prop2.7} claim that in
our circumstances the sets $U_p$ and $U_P$ have the same Hausdorff dimensions
as if we were prescribing frequencies not of the whole $\ell$-words or pairs
of such words but just of their first digits or pairs of their first digits.
\end{remark}

\begin{remark}\label{rem2.9}
It is easy to see that unless $\sum_\be p_{\al\be}=\sum_{\be}p_{\be\al}$
the set $U_P$ is empty, and so the condition (ii) in Proposition \ref{prop2.7}
is a necessary one.
\end{remark}

Next, we consider a similar to Proposition \ref{prop2.5} problem concerning
integer digits $a_0(x),a_1(x),...>0$ of infinite continued fraction expansions
\[
\cfrac{1}{a_0(x)+\cfrac{1}{a_1(x)+\cfrac{1}{a_2(x)+...}}}
\]
for irrational numbers $x\in(0,1)$. We define again $N_\al(x,n)$ and $U_p$
by (\ref{2.21}) and (\ref{2.22}) taking into account that now there are
infinitely many words $\al=(\al_1,...,\al_\ell)\in\{ 1,2,3,...\}^\ell=
\cA_\ell$ and, correspondingly, we have to prescribe infinitely many 
frequencies $p_\al\geq 0,\,\al\in\cA_\ell$ with $\sum_{\al\in\cA_\ell}p_\al
=1$. We recall that the Gauss map $Tx=\{\frac 1x\}$ acts so that $a_i(Tx)
=a_{i+1}(x),\, i=0,1,2,...$, and so $N_\al(x,n)$ can be represented again
in the form (\ref{2.23}). For each infinite probability vector 
$\bar r=(r_1,r_2,...)$ denote by $\cN(\bar r)$ the set of $T$-invariant 
ergodic probability measures $\mu$ such that
\begin{equation}\label{2.28}
\int|\log x|s\mu(x)<\infty\,\,\mbox{and}\,\,\mu[(j+1)^{-1},j^{-1})=r_j\,\,
\mbox{for all}\,\,j\geq 1.
\end{equation}
Here, $[(j+1)^{-1},j^{-1})=\{ x\in(0,1):\, a_0(x)=j\}$ and for any $n$ we
set $I(i_0,i_1,...,i_{n-1})=\{ x\in(0,1):\, a_0(x)=i_0,...,a_{n-1}(x)=
i_{n-1}\}$ which is called a rank-$n$ basic interval. Denote by 
$\hat\cN(\bar r)$ the subset of $\cN(\bar r)$ consisting of measures $\nu$ 
such that for $\nu$-almost all $x$ and all $\al=(\al_1,...,\al_\ell)
\in\cA_\ell$,
\begin{equation}\label{2.29}
\lim_{n\to\infty}\frac 1nN_\al(x,n)=\prod_{i=1}^\ell r_{\al_i}.
\end{equation}
By (\ref{2.23}) and (\ref{2.28}) we see that $\hat\cN(\bar r)$ contains
all measures $\nu\in\cN(\bar r)$ with sufficient mixing which make the
process $X_\al(n)=X_\al(x,n)=(\bbI_{\Gam_{\al_1}}(T^nx),...,
\bbI_{\Gam_{\al_\ell}}(T^nx))$ on the probability space $((0,1),\nu)$ to
 satisfy conditions of Theorem \ref{thm2.2}. We observe that not only the
 Gauss measure $G(\Gam)=\frac 1{\ln 2}\int_\Gam\frac {dx}{1+x}$, which 
 is exponentially fast $\psi$-mixing according to \cite{He}, but also 
 many other $T$-invariant Gibbs measures constructed in \cite{Wa} 
 have sufficiently good mixing properties to satisfy conditions of
 Theorem \ref{thm2.2}. Actually, the rank-1 basic intervals form a Markov
 partition for $T$ whose action is essentially equivalent to the full
 shift on a sequence space with infinite alphabet. For such Markov
 transformations Corollary 4.7.8 from \cite{Aa} gives conditions for
 their Gibbs invariant measures to be exponentially fast $\psi$-mixing.
 
 \begin{proposition}\label{prop2.9}
 Suppose that $q_1(k)=k$ and there exists an infinite probability vector
 $\bar r=(r_0,r_1,...)$ such that $p_\al=\prod^\ell_{i=1}r_{\al_i}$ for
 any $\al=(\al_1,...,\al_\ell)\in\cA_\ell$. Then
 \begin{equation}\label{2.30}
 \max\big(\frac 12,\,\sup_{\nu\in\hat\cN(\bar r)}\frac {h_\nu}
 {2\int|\ln x|d\nu(x)}\big)\leq HD(U_p)\leq\max\big(\frac 12,\,
 \sup_{\nu\in\cN(\bar r)}\frac {h_\nu}{2\int|\ln x|d\nu(x)}\big)
 \end{equation}
 where $h_\mu$ denotes the entropy of $T$ with respect to $\mu$ and "$\sup$"
 is set to be zero if $\cN(\bar r)=\emptyset$.
 \end{proposition}

 \begin{remark}\label{rem2.11}
 All results of this paper can be extended under appropriate conditions 
 to random transformations and processes in random (dynamical) environment.
 Namely, suitable (random) mixing conditions can be introduced similarly to
 \cite{Ki2} and the corresponding relative strong law of large numbers can
 be proved relying on martingale approximations constructed combining methods
 of \cite{Ki2} and \cite{KV}. A relative version of Proposition \ref{prop2.5} 
 can be proved in the spirit of random base expansions from \cite{Ki1}.
 \end{remark}

 \section{Mixingale representation and proof of SLLN}\label{sec3}
\setcounter{equation}{0}

We rely on the following result which is part of Corollary 3.6 from \cite{KV}. 

\begin{lemma}\label{lem3.1}
Let $\cG$ and $\cH$ be $\sig$-subalgebras on a probability space $(\Om,\cF,P)$,
$X$ and $Y$ be $d$-dimensional random vectors and $f=f(x,\om),\, x\in\bbR^d$ be a 
collection of random variables measurable with respect to $\cH$ and satisfying
\begin{equation}\label{3.1}
\|f( x,\omega)-f( y,\omega)\|_{q}\le C (1+|x|^\iota + |y|^\iota)|x-y|^\ka
\,\,\mbox{and}\,\,\|f(x,\omega)\|_{q}\le C (1+|x|^\iota)
\end{equation}
where  $g\geq 1$. 
Set $g(x)=Ef(x,\om)$. Then
\begin{equation}\label{3.2}
\| E(f(X,\cdot)|\cG)-g(X)\|_\up\leq c(1+\| X\|^{\io+2}
_{b(\io+2)})(\vp_{q,p}(\cG,\cH)+\| X-E(X|\cG)\|^\te_{q})
\end{equation} 
provided $\ka-\frac dp>\te>0$, $\frac 1\up\geq\frac 1p+
\frac 1{b}+\frac {\te}q$ with $c=c(C,\io,\ka,\te,p,q,\up,d)>0$ 
depending only on parameters in brackets. Moreover, let
$x=(v,z)$ and $X=(V,Z)$, where $V$ and $Z$ are $d_1$ and $d-d_1$-dimensional
random vectors, respectively, and let $f(x,\om)=f(v,z,\om)$ satisfy (\ref{3.1})
in $x=(v,z)$. Set $\tilde g(v)=Ef(v,Z(\om),\om)$. Then
\begin{eqnarray}\label{3.3}
&\| E(f(V,Z,\cdot)|\cG)-\tilde g(V)\|_\up\leq c(1+\| X\|^{\io+2}
_{b(\io+2)})\\
&\times\big(\vp_{q,p}(\cG,\cH)
+\| V-E(V|\cG)\|^\te_{q}+\| Z-E(Z|\cH)\|^\te_{q}\big).\nonumber
\end{eqnarray}
\end{lemma}

Set $\bar Y_i(n)=F_i(X(q_1(n)),...,X(q_i(n)))-EF_i(X(q_1(n)),...,X(q_i(n)))$
and denote $\cG^{(i)}_n=\cF_{-\infty,q_i(n)}$ for $n\geq 0$ while taking
$\cG^{(i)}_n$ to be the trivial $\sig$-algebra $\{\emptyset,\Om\}$ for
$n<0$. Then by (\ref{2.5}), (\ref{2.6}), (\ref{2.15}) and (\ref{3.3}) of
Lemma \ref{lem3.1} we obtain that for some $C_1>0$ and all $n,m$ and 
$i=1,...,\ell$,
\begin{equation}\label{3.4}
\|E(\bar Y_i(n)|\cG^{(i)}_{n-m})\|_2\leq C_1\big(\vp_{p,q}(\rho_i(m,n))+
\be^\te_q(\rho_i(m,n))\big)
\end{equation}
where $p,q,\te$ satisfy conditions of Assumption \ref{ass2.1} and
\[
\rho_i(m,n)=\min\big(\big[\frac {q_i(n)-q_{i-1}(n)}3,\big]\,
\big[\frac {q_i(n)-q_i(n-m)}3\big]\big).
\]
Observe that $E(\bar Y_i(n)|\cG^{(i)}_{n-m})=0$ if $m>n$ and if $0\leq m\leq n$ 
then $\rho_i(m,n)\geq [\ve m/3]$ by (\ref{2.7}). Hence,
\begin{equation}\label{3.5}
\|E(\bar Y_i(n)|\cG^{(i)}_{n-m})\|_2\leq C_1\big(\vp_{p,q}([\ve m/3])+
\be^\te_q([\ve m/3])\big).
\end{equation}
It follows also from (\ref{2.3}) and (\ref{2.5})--(\ref{2.7}) and the
H\" older inequality (see Lemmas 4.1 and 4.2 together with Theorem 4.4
from \cite{KV}) that
\begin{eqnarray}\label{3.6}
&\big\|\bar Y_i(n)-E(\bar Y_i(n)|\cG_{n+m})\|_2\leq K\big\| \big(1+
\sum_{i=1}^\ell(|X(q_i(n))|^\io\\
&+|E(X(q_i(n))|\cG^{(i)}_{n+m})|^\io)\big)\sum_{i=1}^\ell\big\vert 
X(q_i(n))-E(X(q_i(n))|\cG^{(i)}_{n+m})\big\vert^\ka\big\|_2\leq 
C_2\be_q^\te(m)\nonumber
\end{eqnarray}
 for some $C_2>0$. The estimates (\ref{3.5}) and (\ref{3.6}) yield that
 $\bar Y_i(n),\, n\geq 1$ is a mixingale sequence as defined in \cite{ML}
 and under Assumption \ref{ass2.1} the conditions of Corollary 1.9 from 
 there are satisfied yielding that with probability one for $i=1,...,\ell$,
 \begin{equation}\label{3.7}
 \lim_{N\to\infty}\frac 1N(\Xi_i(N)-E\Xi_i(N))=0.
 \end{equation}
 Set $a_i(n)=(q_{i-1}(n)+q_i(n))/2$. By (\ref{2.11}) and (\ref{3.3}) we
 obtain that
 \begin{eqnarray}\label{3.8}
 &\quad\big\vert EF_i(X(q_1(n)),...,X(q_i(n)))\big\vert=\big\vert EE\big(
 F_i(X(q_1(n)),...,X(q_i(n)))|\cF_{-\infty,a_i(n)}\big)\big\vert\\
 &\leq C\big(\vp_{q,p}([\frac 14(q_i(n)-q_{i-1}(n))])+\be^\del_q([\frac
  14(q_i(n)-q_{i-1}(n))])\big)\to 0\,\,\mbox{as}\,\, n\to\infty\nonumber
  \end{eqnarray}
  by (\ref{2.7}) and Assumption \ref{ass2.1}. It follows that for 
  $i=1,...,\ell$,
  \begin{equation}\label{3.9}
  \lim_{N\to\infty}\frac 1NE\Xi_i(N)=0,\,\,\mbox{and so}\,\,
  \lim_{N\to\infty}\frac 1N\Xi_i(N)=0\,\,\mbox{a.s.}
  \end{equation}
  by (\ref{3.7}) yielding Theorem \ref{thm2.2} in view of (\ref{2.12}).  \qed
 
 In order to derive Corollary \ref{cor2.3} we recall that H\" older 
 continuous functions can be uniformly approximated by functions which
 are constant on elements $\cap_{i=-n}^nT^iG_{k_i}$ of the partition 
 $\bigvee_{i=-n}^nT^i\zeta$ (where $G_k$ are elements of a Markov partition
 $\zeta$) with an error decaying exponentially fast in $n$. Thus
 Theorem \ref{thm2.2} holds when $X_j(n)=f_j(T^nx)$ and $f_j,\, j=1,...,\ell$
 are H\" older continuous. Then Theorem \ref{thm2.2} holds true also for
 continuous functions $f_j,\, j=1,...,\ell$ since they can be uniformly 
 approximated by H\" older continuous ones and $F$ is H\" older continuous.
 Next, let $f_j=\bbI_{\Gam_j}$ with $\mu(\partial\Gam_j)=0$. Given a Markov
 partition $\zeta$ denote by $\tilde\Gam_j^{(l)}$ the set consisting of
 elements of the partition $\bigvee_{i=-l}^lT^i\zeta$ which intersect $\Gam_j$.
 Here we assume that $\Gam_j$ lie on a hyperbolic invariant set itself
 though the argument can be easily extended to a neighborhood of a hyperbolic 
 attractor. For any $\del>0$ there exists $l_\del$ such that 
 $\mu(\tilde\Gam_j^{(l)}\setminus\Gam_j)<\del$ for each $l\geq l_\del$.
 For such an $l$ set $g_j=\bbI_{\tilde\Gam_j^{(l)}},\, j=1,...,\ell$. 
 Since $F$ is H\" older continuous we obtain that
 \begin{eqnarray}\label{3.10}
 &\big\vert\sum_{n=1}^NF(f_1(T^nx),...,f_\ell(T^nx))-
 \sum_{n=1}^NF(g_1(T^nx),...,g_\ell(T^nx))\big\vert \\
 &\leq C\sum_{j=1}^\ell\sum_{n=1}^N(\bbI_{\tilde\Gam_j^{(l)}}(T^nx)-
 \bbI_{\Gam_j}(T^nx))\nonumber
 \end{eqnarray}
 and the remaining part of Corollary \ref{cor2.3} follows by the ergodic 
 theorem applied to the right hand side of (\ref{3.10}). \qed
 
 In the continuous time case we set
 \begin{equation}\label{3.11}
 \bar Y_i(n)=\int_n^{n+1}\big( F_i(X(q_1(s)),...,X(q_\ell(s)))-
 EF_i(X(q_1(s)),...,X(q_\ell(s)))\big)ds
 \end{equation}
 and check similarly to (\ref{3.4})--(\ref{3.6}) that $\big(\bar Y_i(n),
 \cG_n^{(i)}\big)^\infty_{n=1}$ is a mixingale sequence where $\cG_n^{(i)}$
 is the same as in (\ref{3.4})--(\ref{3.6}). Now Corollary \ref{cor2.4}
 follows from \cite{ML} as before. \qed
 
 \section{Application to fractal dimensions}\label{sec4}
\setcounter{equation}{0}

Set 
\[
N_i^{(1)}(x,n)=\#\{ k>0,k\leq n:\, a_k(x)=i\},
\]
\[
N_{ij}^{(1)}(x,n)=\#\{ k>0,k\leq n:\, a_k(x)=i,\, a_{k+1}(x)=j\},
\]
\[
U^{(1)}_r(x,n)=\{ x\in(0,1):\,\lim_{n\to\infty}\frac 1nN_i^{(1)}(x,n)=r_i\,\,
\mbox{for all}\,\, i\}\,\,\,\mbox{and}
\]
\[
U^{(1)}_R(x,n)=\{ x\in(0,1):\,\lim_{n\to\infty}\frac 1nN_{ij}^{(1)}(x,n)
=r_{ij}\,\,\mbox{for all}\,\, i,j\}.
\]
Since 
\[
N_i^{(1)}(x,n)=\sum_{\al_2,...,\al_n}N_{i\al_2...\al_n}(x,n)\,\,\,\mbox{and}
\]
\[
N_{ij}^{(1)}(x,n)=\sum_{\al_2,...,\al_n,\be_2,...,\be_n}N_{i\al_2...\al_n;
j\be_2...\be_n}(x,n)
\]
then for 
\[
r_i=\sum_{\al_2,...,\al_n}p_{i\al_2...\al_n}\,\,\mbox{and}\,\, r_{ij}=
\sum_{\al_2,...,\al_n,\be_2,...,\be_n}p_{i\al_2...\al_n;j\be_2...\be_n}
\]
we obtain that
\[
U_p(x,n)\subset U_r^{(1)}(x,n)\,\,\mbox{and}\,\, U_P(x,n)\subset U_R^{(1)}(x,n)
\]
provided $p=(p_\al,\,\al\in\cA_\ell)$ and $P=(p_{\al\be},\,\al,\be\in\cA_\ell)$.
 Hence, the upper bounds of Propositions \ref{prop2.5}, 
\ref{prop2.7} and \ref{prop2.9} follow from the corresponding upper bounds
from \cite{Bi1}, \cite{Eg} and \cite{FLM}. Still, we provide below an argument
yielding the upper bounds in Propositions \ref{prop2.5} and \ref{prop2.7} by
the reason explained in Remark \ref{rem4.1}.

Denote by $\Xi$ the space of sequences $\xi=(\xi_0,\xi_1,...)$ with $\xi_i
\in\{ 0,1,...,m-1\}$ for all $i\geq 0$. For each probability vector 
$r=(r_0,r_1,...,r_{m-1})$ denote by $\mu_r=(r_0,r_1,...,r_{m-1})^{\bbN}$
the corresponding product measure on $\Xi$, i.e. the probability mesure 
which gives the weight $r_{\al_0}r_{\al_1}\cdots r_{\al_n}$ to each cylinder
set $\Xi_{\al_0\al_1...\al_n}=\{\xi=(\xi_0,\xi_1,...)\in\Xi:\,\xi_i=\al_i\,\,
\mbox{for}\,\, i=0,1,...,n\}$. Observe that the map $\vf:\Xi\to[0,1]$ acting
by the formula $\vf(\xi)=\sum^\infty_{i=1}m^{-i}\xi_{i-1}$ is one-to-one 
except for a countable set of points and since $\mu_r$ has no atoms $\vf$ maps
$\mu_r$ to an atomless measure $\vf\mu_r$ on $[0,1]$. Since $\mu_r$ is invariant
with respect to the left shift $\te:\Xi\to\Xi$ acting by $\te(\xi)=\tilde\xi$
with $\tilde\xi_i=\xi_{i+1}$ then $\vf\mu_r$ is invariant with respect to
$Tx=\{ mx\}$ and $\vf$ provides an isomorphism between $(\Xi,\mu_r,\te)$ and
$([0,1],\vf\mu_r,T)$. Clearly, the conditions of Theorem \ref{thm2.2} are
satisfied here and applying it (see also Remark \ref{rem2.6}) we conclude from
(\ref{2.23}) that for $\vf\mu_r$ almost all $x\in[0,1]$,
\begin{eqnarray}\label{4.1}
 &\lim_{n\to\infty}\frac 1nN_\al(x,n)=\lim_{n\to\infty}\frac 1n\sum_{k=1}^n
 \bbI_{\Gam_{\al_1}}(T^{q_1(k)}x)\bbI_{\Gam_{\al_2}}(T^{q_2(k)}x)\\
 &\times\cdots\times\bbI_{\Gam_{\al_\ell}}(T^{q_\ell(k)}x)
 =\prod_{i=1}^\ell\vf\mu_r(\Gam_{\al_i})=\prod_{i=1}^\ell r_{\al_i}.\nonumber
 \end{eqnarray}
 It follows that
 \begin{equation}\label{4.2}
 \vf\mu_r(U_{p})=1\quad\mbox{if}\quad p=(p_\al,\,\al\in\cA_\ell)\quad\mbox{and}
 \quad p_\al=\prod_{i=1}^\ell r_{\al_i}\,\,\mbox{whenever}\,\,\al=(\al_1,...,
 \al_\ell).
 \end{equation}
 
 Suppose that $r_{i_j}>0,\, j=1,...,k$ while $r_i=0$ if $i\ne i_j$ for any $j$.
 Set
 \[
 U^+=\{ x\in U_p:\, a_j(x)\in\{ i_1,i_2,...,i_k\}\,\,\mbox{for any}\,\, 
 j=0,1,2...\}.
 \]
Then by (\ref{4.2}) and the definition of $\mu_r$,
\begin{equation}\label{4.3}
\vf\mu_r(U^+)=1.
\end{equation} 
Observe that 
 $I_{\al_0\al_1...\al_{n-1}}=\vf \Xi_{\al_0\al_1...\al_{n-1}}$ is a subinterval
 of $[0,1]$ and let $I_n(x)=I_{\al_0\al_1...\al_{n-1}}$ if $x\in 
 I_{\al_0\al_1...\al_{n-1}}$. Set $m_j(x,n)=\#\{ i\geq 0,i<n:\, a_i(x)=j\}$.
 Then for any $x\in U^+$ we can write
 \begin{equation}\label{4.4}
 \ln\vf\mu_r(I_n(x))=\sum_{j=0}^{m-1}m_j(x,n)\ln r_j.
 \end{equation}
 Clearly, $|I_n(x)|=m^{-n}$ where $|I|$ denotes the length of $I$. Observe
 that if $x\in U_p$ with $p=(r_0,r_1,...,r_{m-1})$ then 
 \begin{equation}\label{4.5}
 \lim_{n\to\infty}\frac 1nm_j(x,n)=\lim_{n\to\infty}\frac 1n
 \sum_{0\leq\al_2,...,\al_\ell\leq m-1}N_{j\al_2,...,\al_\ell}(x,n)=r_j.
 \end{equation}
 Hence, for any $x\in U^+$,
 \begin{equation}\label{4.6}
 \lim_{n\to\infty}\frac {\ln\vf\mu_r(I_n(x))}{\ln |I_n(x)|}=-\frac 
 {\sum_{j=0}^{m-1}r_j\ln r_j}{\ln m}
 \end{equation}
 which together with (\ref{4.2}) implies (see Theorem 14.1 in
 \cite{Bi2} or Section 10.1 in \cite{Fa}) that
 \begin{equation}\label{4.7}
 HD(U_p)\geq HD(U^+)=-\frac {\sum_{j=1}^kr_{i_j}\ln r_{i_j}}{\ln m}=
 -\frac {\sum_{i=1}^{m-1}r_i\ln r_i}{\ln m}
 \end{equation}
 with the convention $0\ln 0=0$.
 
 Set $l=m-k$ which is the number of $j\in\{ 0,1,...,m-1\}$ such that
 $r_j=0$. Choose $\del>0$ so small that
 \begin{equation}\label{4.8}
 r_j>\del k^{-1}\,\,\mbox{if}\,\, r_j>0\,\,\ln(\del l^{-1})\leq k^{-1}
 \sum_{j:r_j>0}\ln(r_j-\del k^{-1}).
 \end{equation}
 Set $r^{(\del)}=(r^{(\del)}_0,r^{(\del)}_1,...,r^{(\del)}_{m-1})$ where
 $r^{(\del)}_j=r_j-\del k^{-1}$ if $r_j>0$ and $r_j^{(\del)}=\del l^{-1}$
 if $r_j=0$. Observe that by (\ref{4.8}),
 \begin{equation}\label{4.9}
 \sum_{j=0}^{m-1}r_j\ln r_j^{(\del)}\geq\sum_{j=0}^{m-1}r^{(\del)}_j
 \ln r^{(\del)}_j.
 \end{equation}
 Set
 \[
 W^{(\del)}=\big\{ x\in[0,1]:\,\limsup_{n\to\infty}\big(-\frac 1n\ln\vf
\mu_{r^{(\del)}}(I_n(x))\big)\leq -\sum_{j=0}^{m-1}r^{(\del)}_j\ln r^{(\del)}_j
\big\}.
 \]
 where $\mu_{r^{(\del)}}$ is the Bernoulli measure constructed by $r^{(\del)}$
 in the same way as $\mu_r$ is constructed by $r$. As in (\ref{4.4}),
 \begin{equation}\label{4.10}
 \ln\vf\mu_{r^{(\del)}}(I_n(x))=\sum_{j=0}^{m-1}m_j(x,n)\ln r_j^{(\del)},
 \end{equation}
 and so by (\ref{4.4}), (\ref{4.5}), (\ref{4.9}) and (\ref{4.10}),
 \begin{equation}\label{4.11}
 U_p\subset\{ x\in[0,1]:\,\lim_{n\to\infty}\big(-\frac 1n\ln\vf
 \mu_{r^{(\del)}}(I_n(x))\big)=-\sum_{j=0}^{m-1}r_j\ln r_j^{(\del)}\}\subset
 W^{(\del)}.
 \end{equation}
 If $U_{p^{(\del)}}$ is constructed by $p^{(\del)}=(p^{(\del)}_\al,\,\al\in
 \cA_\ell)$ with $p^{(\del)}_\al=\prod_{i=1}^\ell r^{(\del)}_{\al_i}$ and
 $\al=(\al_1,...,\al_\ell)$ in the same way as $U_p$ is constructed by $p$
 then similarly to (\ref{4.2}) it follows that $\vf\mu_{r^{(\del)}}
 (U_{p^{(\del)}})=1$ and since $U_{p^{(\del)}}\subset W^{(\del)}$ we
 conclude from here and (\ref{4.11}) that
 \begin{equation}\label{4.12}
 \vf\mu_{r^{(\del)}}(W^{(\del)})=1\quad\mbox{and}\quad HD(U_p)\leq 
 HD(W^{(\del)}).
 \end{equation}
 Again, since $|I_n(x)|=m^{-n}$ then it follows from the definition of
 $W^{(\del)}$ by the well known argument (see Theorem 2.3 in \cite{Bi1}
 or the proof of Theorem 14.1 in \cite{Bi2} or Proposition 4.9 in \cite{Fa}
 which also can be adapted to our situation) that 
 \begin{equation}\label{4.13}
 HD(W^{(\del)})\leq -\frac {\sum_{j=0}^{m-1}r_j^{(\del)}\ln r_j^{(\del)}}
 {\ln m}.
 \end{equation}
 Letting $\del\to 0$ we obtain
 \[
 HD(U_p)\leq -\frac {\sum_{j=0}^{m-1}r_j\ln r_j}{\ln m}
 \]
 which together with (\ref{4.7}) completes the proof of Proposition 
 \ref{prop2.5}.
 \qed

 \begin{remark}\label{rem4.1} Many papers and several books disregard the 
 fact that the argument in the first part of the proof above due to Billingsley
 works
 only when all $r_j$'s are positive while without this assumption it leads only
 to the lower bound of dimension. This gap was noticed and repaired first only
 in \cite{Ki1} (though it appeared in later papers, as well). The problem here 
 is that when,
 say, $r_{j_0}=0$ then $\vf\mu_r(I_n(x))=0$ provided $a_i(x)=j_0$ for some 
 $i<n$ and for such $x$ the right hand side of (\ref{4.4}) becomes $-\infty$
 which leads nowhere. In other words, the measure $\vf\mu_r$ "disregards"
 such points while, on the other hand, the set of points $x$ which have zero 
 frequency of appearences of $j_0$ in their $m$-expansions is not countable
 and it cannot be disregarded in the Hausdorff dimension computation. In order
 to prove the result for general probability vectors $(r_0,...,r_{m-1})$
 it is necessary to obtain here an appropriate upper bound for the Hausdorff
 dimension either by a combinatorial argument not related to Billingsley's
 ergodic theory one as in \cite{Eg} or by a simpler perturbation argument 
 above due to my student Z.Hellman which appeared in a more general form in
 \cite{Ki1}.
 \end{remark}
 
 Next, we prove Proposition \ref{prop2.7}. Since for some $n$ the matrix $R^n$
 is a positive matrix then, clearly, each $q_i=\sum_jr_{ij}$ must be positive 
 and for each $i,j$ there exists a sequence $i_1,i_2,...,i_{n-1}$ such that
 $r_{ii_1}r_{i_1i_2}\cdots r_{i_{n-1}j}>0$. Then $q_{ii_1}q_{i_1i_2}\cdots
  q_{i_{n-1}j}>0$, and so $Q^n$ is a positive matrix, as well. Clearly,
  $\sum_iq_iq_{ij}=q_j$, and so $q$ is the unique stationary vector of $Q$.
 Set $\Xi_Q=\{\xi=(\xi_0,\xi_1,...):\, q_{i,i+1}>0\,\,\mbox{for all}\,\, 
 i\geq 0\}$ where $Q=(q_{ij},\, i,j=0,1,...,m-1)$.
 Let $\mu_Q$ be the Markov measure on $\Xi_{q,Q}$ which assigns the weight
 $q_{\al_0}q_{\al_0\al_1}q_{\al_1\al_2}\cdots q_{\al_{n-1}\al_n}$ to each
 cylinder set $R_{\al_0\al_1...\al_n}$ with all $\al_i\in\cA_+$. Then $\mu_Q$
 is invariant with respect to the left shift on $\Xi_{q,Q}$ and its image
 $\vf\mu_Q$ on $[0,1]$ is invariant with respect to $T$. Under assumptions
 of Proposition \ref{prop2.7} the probability matrix $Q$ is a transition matrix
 of an exponentially fast $\psi$-mixing (finite) Markov chain (satisfying
 Doeblin's condition), and so the conditions of Theorem \ref{thm2.2} hold 
 true here. We can also rely on Corollary \ref{cor2.3} since $\mu_Q$ is a
 Gibbs measure for the left shift on $\Xi$ constructed by the function
 $\psi(\xi)=-\ln q_{\xi_0\xi_1},\,\xi=(\xi_0,\xi_1,...)$ (see \cite{Bo}).
 Since $\vf\mu_Q(\Gam_{ij})=q_iq_{ij}=r_{ij}$ we conclude from here together
 with (\ref{2.16}), (\ref{2.26}) and the definition of $U_P$ that $\vf\mu_Q
 (U_P)=1$. If $V_Q=\vf\Xi_Q$ then taking into account that $\mu_Q(\Xi_Q)=1$
 we obtain also that $\vf\mu_Q(U_P\cap V_Q)=1$.
 
 Now, for any $x\in V_Q$ and $I_n(x)$ as above
 \begin{equation}\label{4.14}
 \ln\vf\mu_Q(I_n(x))=\ln q_{a_0(x)}+\sum_{i,j=0}^{m-1}m_{ij}(x,n)\ln q_{ij}
 \end{equation}
 where $m_{ij}(x,n)=\#\{ k\geq 0,k<n:\, a_{k-1}(x)=i\,\,\mbox{and}\,\,
 a_k(x)=j\}$. If $x\in U_P$ then similarly to (\ref{4.5}),
 \begin{equation}\label{4.15}
 \lim_{n\to\infty}\frac 1nm_{ij}(x,n)=r_{ij}=q_iq_{ij}.
 \end{equation}
 It follows that for any $x\in U_P\cap V_Q$,
 \begin{equation}\label{4.16}
 \lim_{n\to\infty}\frac {\ln\vf\mu_Q(I_n(x))}{\ln |I_n(x)|}=-\frac 
 {\sum^{m-1}_{i,j=0}q_iq_{ij}\ln q_{ij}}{\ln m}.
 \end{equation}
 and so similarly to (\ref{4.7}),
 \begin{equation}\label{4.17}
 HD(U_P)\geq HD(U_P\cap V_Q)=-\frac
 {\sum^{m-1}_{i,j=0}q_iq_{ij}\ln q_{ij}}{\ln m}.
 \end{equation}
 
 For the lower bound above we dealt only with points $x\in V_Q$ where 
 $q_{a_i(x)a_{i+1}(x)}>0$ for all $i\geq 0$. In order to obtain the upper 
 bound we employ
 again a perturbation argument which in this case seems to be new. Let
 $l_i,\, i=0,1,...,m-1$ be the number of $j=0,1,...,m-1$ such that 
 $r_{ij}=0$ and set $k_i=m-l_i$. Choose $\del>0$ so small that for all 
 $i,j=0,1,...,m-1$,
 \begin{equation}\label{4.18} 
 r_{ij}>k_i^{-1}\del\,\,\mbox{if}\,\, r_{ij}>0\,\,\mbox{and}\,\,
  \ln(l_i^{-1}\del)\leq k_i^{-1}\sum_{j:r_{ij}>0}\ln((r_{ij}-k_i^{-1}\del)
  q_i^{-1}).
 \end{equation}
 Set $r^{(\del)}_{ij}=r_{ij}-k_i^{-1}\del$ if $r_{ij}>0$ and $r^{(\del)}_{ij}
 =l_i^{-1}\del$ if $r_{ij}=0$. Observe that $\sum^{m-1}_{j=1}r_{ij}^{(\del)}
 =q_i$ and define $q_{ij}^{(\del)}=r_{ij}^{(\del)}q_i^{-1}$ yielding a positive
 $m\times m$ probability matrix $Q^{(\del)}=( q_{ij}^{(\del)})$. By 
 (\ref{4.18}) we have
 \begin{equation}\label{4.19}
 \sum_{i,j=0}^{m-1}r_{ij}\ln q_{ij}^{(\del)}\geq\sum_{i,j=0}^{m-1}
 r_{ij}^{(\del)}\ln q_{ij}^{(\del)}.
 \end{equation}
 
 Set
 \begin{eqnarray*}
 &W^{(\del)}=\{ x\in[0,1]:\,\limsup_{n\to\infty}\big(-\frac 1n\ln\vf
 \mu_{Q^{(\del)}}(I_n(x))\big)\\
 &\leq -\sum_{i,j=0}^{m-1}
 r_{ij}^{(\del)}\max(1,q_i^{(\del)}q_i^{-1})\ln q_{ij}^{(\del)}
 \end{eqnarray*}
 where $\mu_{Q^{(\del)}}$ is the Markov measure constructed by $Q^{(\del)}$
 and its unique stationary vector $q^{(\del)}$ (i.e. $q^{(\del)}Q^{(\del)}=
 q^{(\del)}$) in the same way as $\mu_Q$ was constructed by $Q$ and $q$.
 As in (\ref{4.14}),
 \begin{equation}\label{4.20}
 \ln\vf\mu_{Q^{(\del)}}(I_n(x))=\ln q^{(\del)}_{a_0(x)}+\sum_{i,j=0}^{m-1}
 m_{ij}(x,n)\ln q_{ij}^{(\del)},
 \end{equation}
 and so by (\ref{4.14}), (\ref{4.15}), (\ref{4.19}) and (\ref{4.20}),
 \begin{equation}\label{4.21}
 U_P\subset\{ x\in[0,1]:\,\lim_{n\to\infty}\big(-\frac 1n\ln\vf\mu_{Q^{\del}}
 (I_n(x))\big)=-\sum_{i,j=0}^{m-1}r_{ij}\ln q_{ij}^{(\del)}\}\subset W^{(\del)}.
 \end{equation}
 Let
 \[
 \hat U=\{ x\in[0,1]:\,\lim_{n\to\infty}\frac 1nN_{\al\be}(x,n)=
 \prod^\ell_{i=1}q_{\al_i}^{(\del)}q_{\al_i\al_j}^{(\del)}\,\,\mbox{for all}
 \,\,\al,\be\in\cA_\ell\}.
 \]
 Then for any $x\in\hat U$,
 \begin{equation}\label{4.22}
 \lim_{n\to\infty}\frac 1nm_{ij}(x,n)=\lim_{n\to\infty}\frac 1n
 \sum_{0\leq\al_2,\be_2,...,\al_\ell,\be_\ell\leq m-1}N_{i,\al_2,...,
 \al_\ell,j,\be_2,...,\be_\ell}(x,n)=q_i^{(\del)}q_{ij}^{(\del)}
 \end{equation}
 and by (\ref{4.20}) for any $x\in\hat U$,
 \begin{equation}\label{4.23}
 \lim_{n\to\infty}\frac 1n\ln\vf\mu_{Q^{(\del)}}(I_n(x))=\sum_{i,j=0}^{m-1}
 q_i^{(\del)}q_{ij}^{(\del)}\ln q_{ij}^{(\del)}=\sum_{i,j=0}^{m-1}
 q_i^{(\del)}q_i^{-1}r_{ij}^{(\del)}\ln q_{ij}^{(\del)}.
 \end{equation}
 By Theorem \ref{thm2.2} (or by Corollary \ref{cor2.3}) we obtain that 
 $\mu_{Q^{(\del)}}(\hat U)=1$ and since $\hat U\subset W^{(\del)}$ by 
 (\ref{4.23}) and the definition of $W^{(\del)}$ it follows that
 $\mu_{Q^{(\del)}}(W^{(\del)})=1$. Relying again on Theorem 2.3
 in \cite{Bi1} (or see the proof of Theorem 14.1 in \cite{Bi2}) we
 conclude that 
 \begin{equation}\label{4.24}
 HD(U_P)\leq HD(W^{(\del)})\leq -\frac 
 {\sum^{m-1}_{i,j=0}r^{(\del)}_{ij}\max(1,q_i^{(\del)}q_i^{-1})
 \ln q^{(\del)}_{ij}}{\ln m}.
 \end{equation}
 Since $q$ is the unique probability vector satisfying $qQ=q$ then 
 $q^{(\del)}\to q$ as $\del\to 0$ and letting $\del\to 0$ in (\ref{4.24})
 we arrive at 
 \begin{equation}\label{4.25}
 HD(U_P)\leq -\frac 
 {\sum^{m-1}_{i,j=0}r_{ij}\ln q_{ij}}{\ln m}.
 \end{equation}
 which together with (\ref{4.17}) completes the proof of Proposition 
 \ref{prop2.7}.
 \qed
 
 Concerning Proposition \ref{prop2.9} we explained already at the beginning
 of this section that the upper bound there follows from the upper bound
 derived in \cite{FLM}. Next, we obtain the lower bound
 \[
 HD(U_p)\geq\frac {h_\nu}{2\int |\ln x|}d\nu(x)
 \]
 for any $\nu\in\hat\cN(\bar r)$ in the same way as in Theorem 1
 from \cite{BH} since in addition to arguments there concerning continued
 fractions themselves we need only that $\nu(U_p)=1$ (actually, already
 $\nu(U_p)>0$ is enough) which follows from (\ref{2.29}). 
 
 The remaining bound $HD(U_p)\geq\frac 12$ can be proved similarly to 
 Section 4 in \cite{FLM}. Namely, we construct first points $z\in U_p$
 with $a_n(z)\leq n$ for all $n\geq 1$. In order to do this choose probability
 vectors $(r_1^{(n)},r_2^{(n)},...)$ such that $r_k^{(n)}>0$ when $1\leq k
 \leq n$, $\sum^n_{k=1}p_k^{(n)}=1$ and $\lim_{n\to\infty}r_k^{(n)}=r_k$ for
 any $k\geq 1$. Consider independent integer valued random variables
 $Y_1,Y_2,...$ such that $P\{ Y_n=k\}=r_k^{(n)}$. Applying Theorem \ref{thm2.2}
 we conclude similarly to (\ref{4.1}) that for any $\ell$-word 
 $\al=(\al_1,...,\al_\ell)\in\cA_\ell$ and $P$-almost all $\om$,
 \begin{equation}\label{4.26}
 \lim_{n\to\infty}\frac 1n\sum_{k=1}^n\bbI_{\al_1}(Y_{q_1(k)}(\om))
 \bbI_{\al_2}(Y_{q_2(k)}(\om))\cdots\bbI_{\al_\ell}(Y_{q_\ell(k)}(\om))=
 \prod_{i=1}^\ell r_{\al_i}.
 \end{equation}
 Now, in order to satisfy our conditions we can take any $z$ whose continued 
 fraction expansion have digits $a_n(z)=Y_n(\om),\, n=1,2,...$ with $\om$ 
 such that (\ref{2.26}) holds true. 
 
 Next, let $0\leq m(k)\leq\ell$ be integers
  such that $k^2+m(k)\ne q_i(k)$ for all $k\geq 1$ and $i=1,...,\ell$.
  For $z\in U_p$ constructed above and $b>1$ define the set
  \begin{eqnarray*}
  &G_z(b)=\{ x\in(0,1):\, a_{k^2+m(k)}(x)\in(b^{k^2},2b^{k^2})\,\,\mbox{and}\\
  & a_n(x)=a_n(z)\,\,\mbox{if}\,\, n\ne k^2+m(k)\,\,\mbox{for some}\,\, k\}.
  \end{eqnarray*}
  Then, clearly, $G_z(b)\subset U_p$. Following \cite{FLM} we construct a 
  measure $\mu$ on $G_z(b)$ setting for each rank-$m$ basic interval
  $I_m(x)$ containing $x$,
  \[
  \mu(I_m(x))=\prod_{k=1}^nb^{-k^2}
  \]
  provided $n^2\leq m<(n+1)^2$. Now, in the same way as in \cite{FLM} we can
  show that for any $\te>0$ there exists $b>1$ such that for all $x\in G_z(b)$,
  \[
  \liminf_{r\to 0}\frac {\ln\mu(x-r,x+r)}{\ln r}\geq \frac 12-\te.
  \]
  It follows (see, for instance, Theorem 2.3 in \cite{Bi1} or Proposition
  4.9 in \cite{Fa}) that $HD(U_p)\geq HD(G_z(b))\geq\frac 12-\te$ and since
  $\te>0$ is arbitrary we obtain the required bound. Finally, we observe
  that if $\sum^\infty_{j=1}r_j\ln j=\infty$ then $HD(U_p)=\frac 12$ which
  follows from the latter lower bound and the upper bound of \cite{FLM}.
  \qed

\bibliography{matz_nonarticles,matz_articles}
\bibliographystyle{alpha}

\end{document}